\documentclass[letterpaper,  journal]{IEEEtran}  % Comment this line out
\pdfoutput=1                                                        % if you need a4paper
\IEEEoverridecommandlockouts                              % This command is only

\usepackage{graphics} % for pdf, bitmapped graphics files
\usepackage{graphicx}
\usepackage{subcaption}
\usepackage{array}
\usepackage{cite}
\usepackage{epsfig} % for postscript graphics files
\usepackage{amsmath} % assumes amsmath package installed
\usepackage{amssymb}  % assumes amsmath package installed
\usepackage{color}
\usepackage{float}
\usepackage[noabbrev]{cleveref}
\usepackage{algorithm}
\usepackage{algorithmic}
\usepackage{url}
\usepackage{calrsfs}
\usepackage{textcomp}
\usepackage{multirow}
\usepackage{booktabs}
\DeclareMathAlphabet{\pazocal}{OMS}{zplm}{m}{n}
\newcolumntype{P}[1]{>{\centering\arraybackslash}p{#1}}

\begin{document}

\title{Vulnerabilities of Power System Operations to Load Forecasting Data Injection Attacks}

\author{Yize Chen, Yushi Tan, Ling Zhang and Baosen Zhang \vspace{-4ex}

	\thanks{Y. Chen, Y. Tan, L. Zhang and B. Zhang are with the Department of Electrical and Computer Engineering at the University of Washington,  emails: \{yizechen,ystan,lzhang18,zhangbao\}@uw.edu. This work was supported by University of Washington Clean Energy Institute.}
}

\maketitle

\vspace{-30pt}
\begin{abstract}
We study the security threats of power system operation brought by a class of data injection attacks upon load forecasting algorithms. In particular, with minimal assumptions on the knowledge and ability of the attacker, we design attack data on input features for load forecasting algorithms in a black-box approach.  System operators can be oblivious of such wrong load forecasts, which lead to uneconomical or even insecure decisions in commitment and dispatch. To our knowledge, this paper is the first attempt to bring up the security issues of load forecasting algorithms, and shows that accurate load forecasting algorithm is not necessarily robust to malicious attacks. More severely, attackers are able to design targeted attacks on system operations strategically with additional topology information. We demonstrate the impact of load forecasting attacks on two IEEE test cases.  We show our attack strategy is able to cause load shedding with high probability under various settings in the 14-bus test case, and also demonstrate system-wide threats in the 118-bus test case with limited local attacks. 

\end{abstract}

\section{Introduction}
\label{sec:intro}
%% !TEX root=main.tex
%% !TEX root=main.tex
% The reliable and efficient operation of power systems calls for a detailed knowledge of the future demands.
Load forecasting plays an important role in the planning and operations of electric grids. As a cornerstone application for utilities and operators, it provides  future load information which is utilized for various decision-making problems such as unit commitment, reserve management, economic dispatch and maintenance scheduling~\cite{gross1987short}. Consequently, the accuracy of forecasted loads directly impacts the cost and reliability of system operations~\cite{hobbs1999analysis}. %\todo{With a growing penetration of new technologies into the demand side, the importance of accurate and robust forecasts will continue to increase for utilities and system operators.}

Because of its fundamental importance, there are always strong incentives to improve short-term forecasting methods, especially under higher penetration of renewables. The driving factors of load variations are heterogeneous, including temperature, weather, temporal and seasonal effects (e.g., weekday vs. weekend) and other socioeconomic factors. Thus, load forecasting algorithms can be regarded as finding a  nonlinear and complex mapping between the (potentially high dimensional) driving factors to the forecasted time series of load values. Over the past decades, a myriad of load forecasting algorithms have been proposed and adopted. See, for example, \cite{gross1987short, de200625} and the references within. Statistical and machine learning techniques, such as support vector regression~\cite{ceperic2013strategy}, ARIMA~\cite{contreras2003arima} and neural networks~\cite{hippert2001neural} have been applied to short term load forecasting and implemented in practice. The recent advances in deep learning and data sciences opened the door to utilizing more input features and deeper model architectures to further improve load forecasting accuracy and provided some of the best performances to date~\cite{kong2017short, quilumba2015using,wang2016electric}. In most of the load forecasting studies, \emph{forecast accuracy} has been regarded as the holy grail for researchers.

\begin{figure}[h]
	\centering
	\includegraphics[width=1.0 \columnwidth]{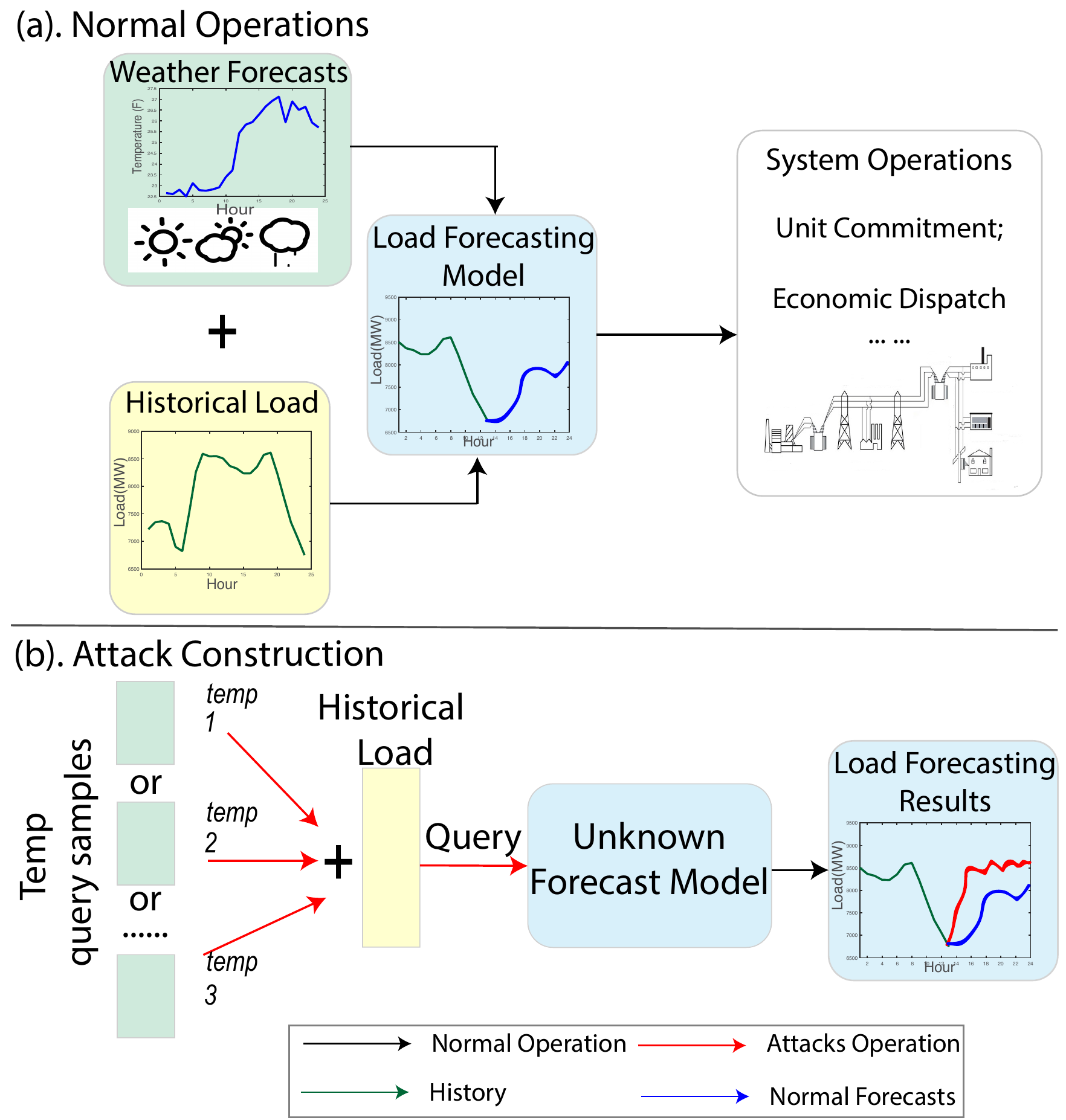}
	\caption{\small The schematic of normal power system operations (a) and proposed black-box attacks on load forecasting (b). Without knowledge about the forecast model's parameters, the attacker injects designed small, undetectable data perturbations into weather forecasts to induce abnormal system operations.} %threats first
	\label{fig:schematic}
\end{figure}

However, despite the critical role of forecasting algorithms and the long pursuit of \emph{forecast accuracy}, the robustness and  security issues have been overlooked to some extents. As the forecasting methods become more complex, they are also more susceptible to cybersecurity threats.  In the previous research on cyber-security of power systems~\cite{mcdaniel2009security, mo2009secure}, where state estimation~\cite{kosut2010malicious, liu2011false}, communication~\cite{ericsson2010cyber}, and electricity market~\cite{xie2010false,tan2018online} threats and countermeasures are rigorously evaluated. However, the vulnerabilities of load forecasting algorithms are rarely discussed~\cite{luo2018benchmarking,chen2019exploiting}, while this does not mean load forecasting is less vulnerable nor the consequences of attacks are  less severe.  For instance, forecasting models normally make use of weather forecasts inputs coming from external services/APIs, while such inputs can be exposed to adversarial modification and the model performance may be severely impacted by such malicious changes. Recently, there has been a hot debate on the security of machine learning models~\cite{szegedy2013intriguing}, and researchers found that small noises injected to the inputs can severely impact model performances~\cite{papernot2016limitations}.
 % In addition, successful distortion on load forecasting also impacts the reliable operations of power systems, thus it is essential to investigate the data vulnerabilities in existing load forecasting methods.

 In this paper, we look into the security threats in general load forecasting algorithms.  By taking the perspective of an attacker and developing attack strategies on load forecasting algorithms, we conduct damage analysis of the proposed attacks. We consider the scenario where attacker adversarially injects false data into the input features of forecasting algorithms, and examine to what extent such attacks could impact the performance of load forecasting models. Specifically, we investigate false data injection attacks on the temperature data. It is an important input to load forecasting algorithms and is mostly obtained from external services/APIs. Therefore it is easier for attackers to inject perturbations on temperature data  than to attack the state estimation~\cite{liu2011false} or market clearing algorithms~\cite{xie2010false}. The potential damage of load forecasting attacks can be significant, leading to increases in system operation costs and maybe even more catastrophic events such as load shedding.

In Figure~\ref{fig:schematic}, we show the schematic of threats and proposed attacks to systems. Our work is different from most related work in two aspects: most of the studies in forecasting lack security and robustness considerations, while  most of the studies in power system security evaluate attacks with certain level of knowledge about the targeted system or unconstrained capabilities. In contrast to previous attacks research that assume full knowledge of system configuration and strong capabilities of attackers on implementing attacks~(see, e.g. \cite{kosut2010malicious, liu2011false,xie2010false}), we take a restrictive setting of both the attacker's \emph{knowledge} and \emph{capabilities}. Under our setup, the attacker does not need to know any parameter of the targeted load forecasting algorithms, and could only inject constrained perturbations into input temperatures to avoid detection.  We develop a simple data-driven attack strategy for finding the injected perturbations onto the original temperature data. Surprisingly, we find the proposed attacks significantly degrade the performance of a class of (accurate) load forecasting algorithms. With only few degrees of perturbations injected into input temperatures, the load forecasting algorithm's output deviates drastically from original accurate forecasts.

To further illustrate security issues brought by load forecasting attacks, we embed the load forecasting algorithms into canonical power system operation case studies considering network constraints and power balance. We consider both cases when attacker possesses and does not possess additional information on system topology and parameters. For the former case where the attacker can strategically inject data perturbations under certain attack budgets, we design a greedy algorithm to compromise a subset of nodal load forecasts to cause targeted damages such as uneconomical generation, infeasible line flows and generator schedules, and load shedding.  Simulations based on real-world load datasets on the  IEEE 14-bus and 118-bus systems demonstrate the system operation vulnerabilities by only maliciously changing the temperature by a few degrees.
%power system operations under adversarial forecasts not only increase the operational costs, but also cause severe damages such as load shedding.

This study illustrates the need to look at other properties in addition to \emph{forecast accuracy}, and the need for more comprehensive analysis when developing and applying load forecasting techniques. We demonstrate that accuracy may not mean robustness, and a wrong forecast of load potentially leads to costly operation decisions or system damage.  Specifically, we make the following contributions in this work:
\begin{itemize}
	\item To the best of our knowledge, this work is the first to evaluate the security issues of load forecasting procedures in power system operations. Starting with the setup for load forecasting along with its role in power system operations (Sec.~\ref{sec:formulation}), data vulnerabilities of current forecasting methods are formulated and discussed  (Sec.~\ref{sec:method}).
	\item Black-box attack algorithm  \texttt{gradient estimation} is proposed to generate hard-to-detect, adversarial input data for load forecasting algorithms  (Sec.~\ref{sec:method}).
	\item We show that the strategically designed adversarial injections upon input features could target either increased system operating costs or load shedding. The resulting optimization formulation of the attack problem maybe of independent interest (Sec.~\ref{sec:market}).
	\item Case studies of power system operations on standard IEEE test cases using real-world load data reveal the prevalent vulnerabilities of current forecasting techniques and demonstrate potential damages on power system operations via proposed attacks (Sec.~\ref{sec:simulation}). 
\end{itemize}

Compared to our prior work~\cite{chen2019exploiting} which works on single bus vulnerabilities analysis, we bring out both the load forecasting threats and the attack strategies on power networks. Extensive numerical simulations also verify such load forecasting vulnerabilities generally exist in power networks operations. We also make our code open source as a public package for evaluating load forecasting robustness and security\footnote{\url{https://github.com/chennnnnyize/load_forecasts_attack}}. Due to the space limits of submission, we refer to the preprint for more details on attack implementations and thorough tests on various load forecasting algorithms \cite{chen2019arxiv}.

%The rest of the paper is organized as follows. We briefly summarize a general load forecasting model, and formulate the objective and constraints of attackers in Section~\ref{sec:formulation}; in Section~\ref{sec:method}, we detail the algorithms for implementing the attack; to illustrate the attack's threats to the power system operations, we describe the market setup and a toy example in Section \ref{sec:market}; through simulations based on real-world load data in Section~\ref{sec:simulation}, we demonstrate the threats posed by the proposed attacks on standard IEEE test system; furthur discussion on model/data security and conclusion are drawn in Section~\ref{sec:conclusion}.

\section{Preliminaries}
\label{sec:formulation}
%% !TEX root=main.tex
In this section, we briefly describe the notations and setup for load forecasting algorithms, and illustrate how load forecasting serves as an important component of system operation in day-ahead commitment and real-time dispatch.

\subsection{Load Forecasting}
To set up and find parameters of the short-term load forecasting algorithm for a specific region, the system operator needs to collect a training dataset $\mathcal{D}_{tr}=\{(\textbf{X}_{t-H},..., \textbf{X}_t);L_{t+k}\}_{t=1}^T$ based on available historical data. Here $L_{t+k}\in [0,1]$ are scalars representing scaled nodal load values and $H$ is the history horizon~\cite{gross1987short}.  The model's forecast horizon is denoted by $k$ and ranges from one hour to one day in short-term forecasts. $\textbf{X}_{t-i}\in [0,1]^d,\; i=0,...,H$ are scaled, $d$-dimensional input feature vectors. Feature vector $\textbf{X}_t$  includes historical records of load, weather forecasts including temperature, weather indicators (e.g., sunny, rainy or cloudy) and seasonal indicator variables such as weekdays/weekends and hour of the day~\cite{wang2016electric}. We express it as $\textbf{X}_t:=\{L_t, \textbf{X}_t^{temp}, \textbf{X}_t^{index}\}$, where $L_t$ is the load history records;  $\textbf{X}_t^{temp}$ is the temperature value vector of current and neighboring regions, which could be acquired from either system historical records or weather forecast API; $\textbf{X}_t^{index}$ are a collection of indicators including weather characteristics, seasonal factors and time factors.  In the task of load forecasting, one is interested to find a function parameterized by $\theta$: $f_\theta(\cdot)$, which learns the mapping from $(\textbf{X}_{t-H},...,\textbf{X}_t)$ to future loads $\hat{L}_{t+k}$.  The mean absolute error~(MAE) is widely used to measure the performance of $f_\theta$, which is defined by the  average $L_1$ norm of difference between forecasted loads $\hat{L}_{t+k}$ and ground truth ${L}_{t+k}$. The schematic of load forecasting model is depicted in Figure~\ref{fig:schematic}(a).

We note that vulnerability analysis conducted by this paper is not constrained to certain forecasting algorithms. As long as the model output is sensitive with respect to input features, our proposed attack methods would be able to alter load patterns maliciously. For the discussion hereafter, we make use of an Recurrent Neural Networks (RNN)~\cite{hippert2001neural,kong2017short, quilumba2015using,wang2016electric}, which is a  widely adopted forecast algorithm to model the temporal dependencies between feature inputs and forecasted values.

\subsection{System Operations Model}
For a $N$-bus power network with set of loads $D$ and set of generators $G$, we consider the system operation setting consisting of a day-ahead planning stage using load forecasts and a real-time operational stage using actual load.

\begin{enumerate}
	\item A unit commitment (UC) model considering reserve margins, startup and shutdown cost,  minimum up/down time constraints and ramping constraints is used to create a commitment schedule $G_t, t=1,...T$ based on the day-ahead load forecasts $\hat{L}_t, t=1,...T$;
	\item For each time $t$, the dispatch of the scheduled units $\mathbf{p}_t$ and the actual dispatch costs $C_1(\mathbf{p}_t)$ are calculated according to a basic Economic Dispatch (ED) model~\cite{kirschen2018fundamentals} based on the actual load $L_t$ and generation schedule $G_t$.
\end{enumerate}

The actual daily operation costs are calculated via summing the 24-hour dispatch costs and the startup and shutdown costs. When ED based on the day-ahead commitment does not have a feasible solution, load shedding is used to maintain the balance between supply and demand. The shedded loads $\mathbf{LS}_t$ also incur costs $C_2(\mathbf{LS}_t)$. Note that under perfect forecasts, the generation schedule $G_t$ can minimize real-time dispatch costs.
 % yet we want to examine if attacks on load forecasts would cause non-optimal or even infeasible schedule for real-time dispatch.

\section{Attack Strategies}
%% !TEX root=main.tex
\label{sec:method}
%% !TEX root=main.tex
In this section, we first describe the objective and constraints for implementing load forecast attacks. We then illustrate how an attacker is able to design \emph{white-box} attack with known load forecasting model parameters.  Finally we describe under the \emph{black-box} setting, how data injection attacks can be implemented even though the attacker only has limited query access to the load forecasting model.

\subsection{Objective of Attacker}
The attacker's goal is to distort the forecasted load as much as possible in a certain direction, e.g., to either increase or decrease forecasted values. Consider the task of training an accurate load forecasting models, where estimation of $\theta$ is given by minimizing the $L_1$-norm of the difference between model predictions and ground truth values:
\begin{equation}
	% \begin{align}
	\min_\theta \frac{1}{T}\sum_{t=1}^{T}||f_{\theta}(\textbf{X}_{t-H},...,\textbf{X}_t)-L_{t+k}||_1 \\
	% s.t.  \quad  & \theta \in \Theta
	% \end{align}
	\label{equ:forecasts}
\end{equation}
where during training, ground truth of historical records on $\textbf{X}_t$ and $L_{t+k}$  are used; during testing and real-world system implementations, we are using $\textbf{X}_t$ which are coming from weather forecast as input features. Once the model is learned, it can be applied in a rolling-horizon fashion.
 % to make use of forecasted $\hat{L}$ along with $\textbf{X}_t^{temp}$ and $\textbf{X}_t^{index}$ to forecast for furthur into the future.

In order to distort the output forecast values from the trained model based on \eqref{equ:forecasts}, the attacker actually has two choices of inserting attacks: \emph{to attack $\textbf{X}_t$} or \emph{to attack $\theta$}. While trained model $\theta$ itself is often safely kept by the operators, system operator has to use external data such as weather forecasts $\textbf{X}_t^{temp}$ as input features for $f_\theta$. This actually provides a backdoor for the attacker, whose goal is to inject perturbations into the weather forecasts coming from external services. By generating adversarial input data $\tilde{\textbf{X}}_t^{temp}$ for $f_\theta(\cdot)$, model predictions are modified adversarially. We use $\gamma=\{-1,1\}$ to denote the chosen attack direction by attackers. If $\gamma=1$ ($\gamma=-1$), the attacker tries to find $\tilde{\textbf{X}}$ to decrease (increase) the load forecasts values.
% when $\gamma=-1$, the attacker tries to find $\tilde{\textbf{X}}$ to  increase load forecasts values.
Since load values are always positive, the attacker's goal is to  find $\tilde{\textbf{X}}$ that minimizes the value of $\gamma f_{\theta}(\tilde{\textbf{X}}_{t-H},...,\tilde{\textbf{X}}_t)$.

\subsection{Attacker's Knowledge}
We consider two attack scenarios, \emph{white box} and \emph{black-box} attacks.
In the \emph{white-box} settings, the attacker is assumed to know exactly the model parameters $\theta$. This is a strong assumption in the sense that load forecast model $f_\theta(\cdot)$ is fully exposed to the attacker. On the contrary, in the \emph{black-box} setting, the attacker only knows which family of load forecasting model has been applied~(e.g., NN or RNN), but is blind to the forecasting algorithms and has no knowledge of any parameters of $f_\theta$. We assume the attacker only has \emph{query} access to the load forecasting model\footnote{Such query access assumption is possible in many \emph{forecast-as-a-Service} businesses, e.g., SAS energy forecasting and Itron forecasting.}.  That is, the attacker could query the implemented load forecasting model by using different values of input features for a limited number of times, and then try to get insights on how $f_\theta$ works.

\subsection{Attacker's Capability}
From the attacker's perspective, it is necessary to construct attack injections while avoid being detected by the system operators' bad data detection algorithms.  We consider several realistic constraints for attacker's capability: it could be upper bounded by the maximum number of perturbed entries in the input data, by the average deviations on all features, or by the largest deviation from the clean data. Mathematically, the attacker wants to keep $\vert\vert\tilde{\textbf{X}}_t^{temp}-\textbf{X}_t^{temp}\vert\vert_p$ bounded, where $p$ can take different values such as $0, 1, \infty$ to express certain norm constraints corresponding to different detection countermeasures.

In summary, we formulate the model of attackers as the following optimization problem:

\begin{subequations}
	\begin{align}
	\label{equ:obj}
	\min_{\tilde{\textbf{X}}_{t-H}^{temp},...,\tilde{\textbf{X}}_t^{temp}}  \quad & \gamma f_{\theta}(\tilde{\textbf{X}}_{t-H},...,\tilde{\textbf{X}}_t)\\
	\label{equ:constraint}
	s.t.  \quad  & ||\textbf{X}_{t-i}^{temp}-\tilde{\textbf{X}}_{t-i}^{temp}||_p \leq \epsilon,\; i=0,...,H
	\end{align}
	\label{equ:attack}
\end{subequations}
%\todo{ forall i? modify not just the current i, but ...}

\begin{figure*}[h]
	\centering
	\includegraphics[width=1.8 \columnwidth]{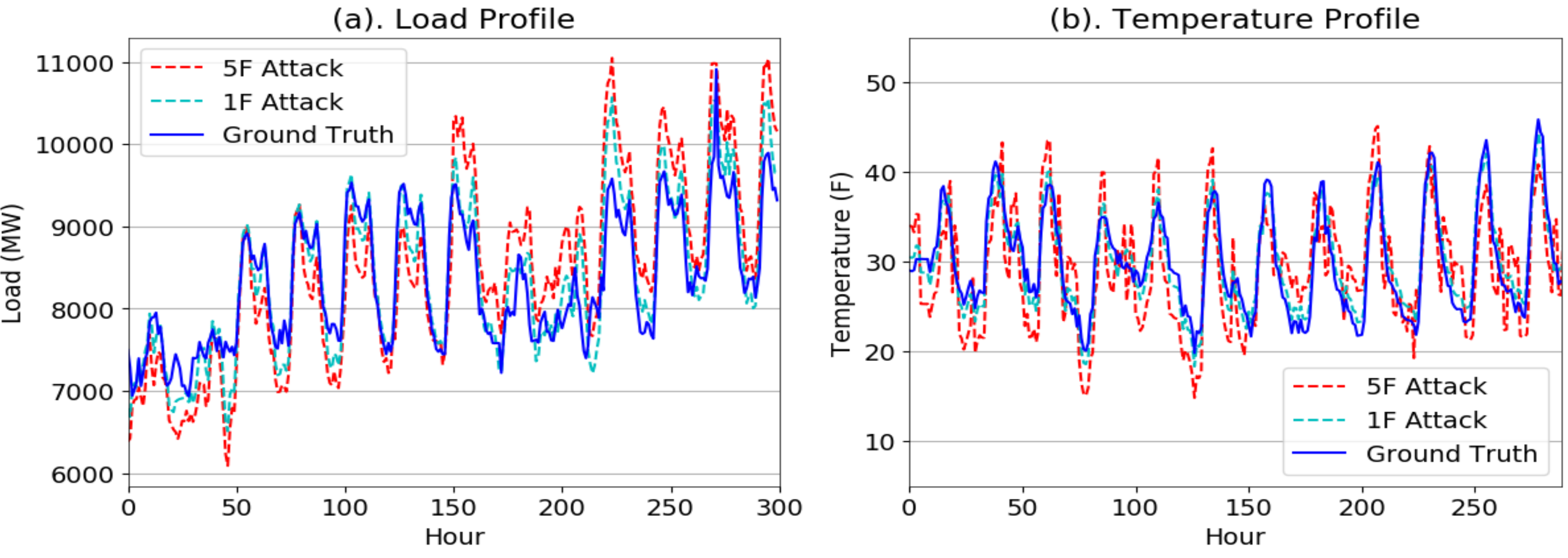}
	\caption{\small Visualization of $300$ hours forecasts based on original and adversarial temperature data for the aggregated load of Switzerland. The load forecasting model is a recurrent neural networks,  while attack perturbations on test dataset are generated by using the \texttt{gradient estimation} method. The attack tries to minimize the load in the first $150$ hours and to maximize the load in the latter $150$ hours. (a). Load forecasting results; (b). false data injections on temperature.}
	\label{fig:forecast_example}
\end{figure*}

Note that there is a parallel between the forecast problem \eqref{equ:forecasts} and attack problem \eqref{equ:attack}, where the objective's optimization directions and optimization variables are exactly in the opposite directions: forecasting model works on model parameters to minimize forecast errors, while attacker works on model inputs to maximize the errors to targeted directions. However, due to lack of model knowledge in the black-box setting, it is a challenging task for attackers to find efficient attack input $\tilde{\textbf{X}}^{temp}$ via \eqref{equ:attack}. In the next subsection, we will show a black-box attack method generally working with attacker's knowledge coming from \emph{query access} to the forecast algorithm.

\subsection{Black-Box Attack}
Under the case of \emph{white-box} where model parameters are known to the attacker, it is possible to find the attack input via solving \eqref{equ:attack}. For the convenience of notations, we omit the superscript on $\textbf{X}$ in some of the following paragraphs, and introduce the generalizable attack methods not only suitable for attacking temperature forecasts, but also suitable for injecting false data into other features.

Since most state-of-the-art load forecasting algorithms use complex models such as neural networks, the attacker's problem \eqref{equ:attack} is nonconvex and furthermore, there is no closed-form solution for $\tilde{\textbf{X}}_{t-H},...,\tilde{\textbf{X}}_t$. Nevertheless, an attacker can still find some attack vectors iteratively by taking gradients with respect to each time step's temperature values. Even though this may not find the optimal solution to \eqref{equ:attack}, because of the highly nonconvex nature of the forecasting model, a slight (suboptimal) perturbation of the input features would drastically change the forecast output.

Based on \eqref{equ:attack}, we define a loss function $\mathcal{L}$ with respect to each time step's feature $\tilde{\textbf{X}}_{t-i}, \; i=0,...,H$. Then the attacker iteratively takes gradients of $\mathcal{L}$ to find the adversarial input $\tilde{\textbf{X}}_{t-i}$. The constraints in \eqref{equ:constraint} is included in the loss function using a log-barrier:
\begin{equation}
\label{equ:loss}
\mathcal{L}(\tilde{\textbf{X}}_{t-i})=\gamma  f_{\theta}(\tilde{\textbf{X}}_{t-H},...,\tilde{\textbf{X}}_t)- \beta \log( \epsilon- ||\textbf{X}_{t-i}^{temp}-\tilde{\textbf{X}}_{t-i}^{temp}||_p)
\end{equation}
where $\beta$ is the weight of the barrier term.
Since there are a large number of parameters and input features in many load forecasting algorithms, it can be computationally expensive to compute the exact gradient values for each input feature. We follow a simpler method in~\cite{szegedy2013intriguing} to only update the attack features based on the sign of the gradient at each iteration $j$:
\begin{equation}
\label{equ:grad_attack}
\tilde{\textbf{X}}_{t-i}^{(j+1)}=\tilde{\textbf{X}}_{t-i}^{(j)}-\alpha \cdot \text{sign}(\triangledown_{\tilde{\textbf{X}}_{t-i}^{(j)}} (\mathcal{L}(\tilde{\textbf{X}}_{t-i}^{(j)})))
\end{equation}
where $\alpha$ controls the step size for updating adversarial temperature values. The resulting adversarial temperature vector is obtained by applying \eqref{equ:grad_attack} for several times.

Under the black-box setting where $\triangledown_{\tilde{\textbf{X}}_{t-i}^{(j)}}(\mathcal{L}(\tilde{\textbf{X}}_{t-i}^{(j)}))$ is not able to compute, we assume attacker is able to query the load forecasting algorithm for a limited number of times, and it is still possible to construct adversarial temperature inputs by using queries to estimate the gradients. In Figure \ref{fig:schematic}(b) we show the schematic on generating adversarial temperature instances via querying. For $k$-th dimension of the input feature at time stamp $t-i$, $\tilde{\textbf{X}}_{k,t-i}^{j+1}$, the attacker needs to query the load forecasting system on each feature entry to calculate the two-sided estimation of the gradient of $f_\theta$:
\begin{equation}
\triangledown_{\tilde{\textbf{X}}_{k,t-i}}f_\theta(\tilde{\textbf{X}})   \approx \frac{f_\theta(\tilde{\textbf{X}}+\delta\mathbf{e}_k)-f_\theta(\tilde{\textbf{X}}-\delta\mathbf{e}_k)}{2\delta}
\end{equation}
where $\mathbf{e}_k$ is a $d$-dimensional vector with all zero except $1$ at $k$-th component, and $\delta$ takes a small value for gradient estimation. Once the gradient is estimated for each dimension of temperature features, we can follow the same method of \eqref{equ:grad_attack} to iteratively build the adversarial features using the estimated gradient vectors:
\begin{equation}
\label{equ:grad_attack2}
\tilde{\textbf{X}}_{t-i}^{(j+1)}=\tilde{\textbf{X}}_{t-i}^{(j)}-\alpha \gamma \cdot \text{sign}(\triangledown_{\tilde{\textbf{X}}_{t-i}}f_\theta(\tilde{\textbf{X}}^{(j)}) ).
\end{equation}

To satisfy norm constraints on the allowed perturbation of $\tilde{\textbf{X}}$, the attacker projects the adversarial data back into the pre-defined norms after each iterative attack step. As shown in Figure~\ref{fig:forecast_example}, the output load forecasts deviate a lot from the ground truth values, while temperature perturbations are constrained within attacker's capability. In~\cite{bhagoji2018practical}, techniques on reducing number of queries are discussed for attacking an image classifier, which may also further improve the query efficiency of load forecasting attacks.

\section{Attacks on System Operations}
%% !TEX root=main.tex
\label{sec:market}
%% !TEX root=main.tex
In this section, we illustrate the attacks on load forecasting input features could further threaten power system operations, and propose a realistic attack strategy for the attacker under attack budgets on number of compromised nodal forecasts.

\subsection{Attack Objectives}
We assume the attacker could only inject constrained attacks in the day-ahead planning stage. Under the day-ahead load forecasts $\hat{L}_d^t$ for $|D|$ loads in the networks, the UC problem is to find the generation schedule and dispatch for $|G|$ generators while satisfying reserve and system operation constraints:

\vspace{-10pt}
\begin{subequations}
	\label{equ:UC}
	\begin{align}\vspace{-20pt}
		\underset{\mathbf{u}, \mathbf{p} }{\min} \;\, & C(\mathbf{p}) + S(\mathbf{u}) \\
		\text{s.t.} & \textstyle \sum_{g \in G} p_{g}^{t} = \sum_{d \in D}\hat{L}_d^t, \;\, \forall t \in T \label{eqn:UC_balance} \\
		& u_i^t p^{\min}_g \leq p_g^t \leq u_g^t p^{\max}_i, \;\, \forall g \in G, \;\, \forall t \in T \label{eqn:UC_GLim} \\
		& f_l^{\min}\leq f_l^t \leq f_l^{\max} , \forall l \in F, \forall t\in T \label{eqn:UC_FLim}\\
		& \textstyle \sum_{g\in g(d)}p_g^t+\sum_{k \in \delta(d)}f_k^t=\hat{L}_d^t , \; \forall d\in D,\; t\in T \label{eqn:UC_Nodal} \\
		& u_g^t - u_g^{t - 1} = z_g^t - y_g^t, \;\, \forall g \in G, \;\, t \in T \label{eqn:UC_GLogic} \\
		& \textstyle \sum_{\tau = t - t_g^{\text{up}} + 1}^t z_g^{\tau} \leq x_g^t, \;\, \forall g \in G, \forall t \in T \label{eqn:UC_minup} \\
		& \textstyle \sum_{\tau = t - t_g^{\text{dn}} + 1}^t z_g^{\tau} \leq 1 - x_g^t, \;\, \forall g \in G, \forall t \in T \label{eqn:UC_mindn} \\
		& -R_g^{\text{dn}} \leq p_g^{t + 1} - p_g^t \leq R_g^{\text{up}}, \;\, \forall g \in G \label{eqn:UC_Ramp} \\
		& u_g^t, z_g^t, y_g^t \in \{0, 1\}, \;\, \forall g \in G, \;\, t \in T
	\end{align}
\end{subequations}
where $\delta(d)$ is the set of lines connected to node $d$; $u_{g}^{t}$ is the binary decision variable of the commitment status of generator $g$ at time $t$, with $1$ indicating $g$ is online; $p_g^t$ is the real power output of generator $g$ at time $t$; all the $u_g^{t}$'s and $p_g^t$'s are collected together into vectors $\mathbf{u}$ and $\mathbf{p}$; $C(\mathbf{p})$ and $S(\textbf{u})$ represent the dispatch costs and startup and shutdown costs, respectively, of all the generators in all periods; the constraints are system-wide power balance constraint~(\ref{eqn:UC_balance}), generation limits constraints~(\ref{eqn:UC_GLim}), line flow limits~(\ref{eqn:UC_FLim}), power balance at each node~(\ref{eqn:UC_Nodal}), generator logical constraint~(\ref{eqn:UC_GLogic}), minimum up time constraint~(\ref{eqn:UC_minup}), minimum down time constraint~(\ref{eqn:UC_mindn}) and ramping constraints~(\ref{eqn:UC_Ramp}). We also keep a fixed reserve margin throughout the simulation. Once solved, the operator gets the schedule for the set of online generators $G_t$ at each time $t$.

The attacker injects $\tilde{\mathbf{X}}$ to a group of compromised load forecasts, such that the nodal load forecasts maliciously change from $\hat{L}_d^t$ to $\tilde{L}_d^t$. The attacker is not only constrained by the average deviations upon $\tilde{\mathbf{X}}$~(Equation \eqref{equ:attack}), but also constrained on number of compromised loads $N_{adv}$ that are allowed to inject perturbations throughout the day. From the attacker's perspective, it is then most harmful to find a constrained adversarial day-ahead commitment schedule $\tilde{G}_t$ via $\tilde{L}_d^t$, which maximizes the solution for the following real-time ED:

\begin{subequations}
		\label{equ:ED}
	\begin{align}
	    \underset{ \mathbf{p}_t, \mathbf{LS}_t  }{\min} \,\; & C(\mathbf{p}_t)+C(\mathbf{LS}_t) \\
		\text{s.t.} & \textstyle \sum_{g \in G_t} p_{g}^{t} + \sum_{d \in D}LS_d^t= \sum_{d \in D}L_d^t, \;\, \label{eqn:ED_balance} \\
		& p^{\min}_g \leq p_g^t \leq p^{\max}_g, \;\, \forall g \in \tilde{G}_t \label{eqn:ED_GLim}\\
		& f_l^{\min}\leq f_l^t \leq f_l^{\max}, \forall l \in F, \forall t\in T \label{eqn:ED_FLim}\\
		& \textstyle \sum_{g\in g(d)}p_g^t+\sum_{k \in \delta(d)}f_k^t= L_d^t , \forall d\in D,t\in T  \label{eqn:ED_Nodal}
	\end{align}
\end{subequations}
where from the system operator's perspective, ED aims to find the real power dispatch at time $t$, $\mathbf{p}_t$, that minimize the dispatch costs at time $t$, $C(\mathbf{p}_t)$, considering system-wide power balance constraint~(\ref{eqn:ED_balance}), generation limits constraints~(\ref{eqn:ED_GLim}), line flow limits~(\ref{eqn:ED_FLim}) and power balance at each node~(\ref{eqn:ED_Nodal}).

\subsection{Attack Strategies}
\label{subsec:strategy}
Under normal operating conditions, the load forecasting algorithms provide accurate forecasts on day-ahead load for system operators to solve~\eqref{equ:UC}. When the system is under attack, the attacker chooses a group of load buses to inject  adversarial temperature forecasts, such that generation schedule coming from  day-ahead planning stage is deviating from the normal schedule. Such adversarial generation schedules are likely to cause malicious operation patterns, e.g., increased system costs, load shedding, no feasible generation dispatch or violation of ramping constraints. Essentially, the attacker wants to answer the following questions to find the attacks:
\begin{itemize}
	\item Which group of load buses should be compromised to inject $\tilde{\mathbf{X}}$?
	\item How to generate $\tilde{L}_j^t,t=1,...,T$ for compromised load bus $j$, such that \eqref{equ:ED} is maximized?
\end{itemize}

Under the case all system parameters are known, the attacker's optimization problem is tri-level with integer constraints, which is a very challenging problem to solve. Rather than the standard approach through KKT conditions, we design a modified greedy search algorithm for the attacker to find the most vulnerable loads that cause system-level misoperations.  As described in Algorithm~\ref{algorithm}, the attacker follows a modified best-first search algorithm to implement attacks on compromised nodal forecasts iteratively~\cite{dechter1985generalized}. For each iteration, the attacker checks the lines and generators which are approaching operating limits (e.g., line flow capacity, generation capacity), and finds the most vulnerable load $j$ based on neighboring lines and generators. Then the attacker raises $\tilde{\underline{\mathbf{X}}}_j$ and $\tilde{\bar{\mathbf{X}}}_j$ which maximizes and minimizes the load forecasts respectively. By solving \eqref{equ:UC} using the resulting  $\tilde{\underline{L}}_j$ and $\tilde{\bar{L}}_j$, the attacker checks which candidate attack is more prone to make day-ahead commitment schedule $\tilde{G}_t$ different from $G_t$ and keeps it as $\tilde{L}_j$ in this iteration, while $j$ is added to the set of compromised loads $\tilde{\mathbb{L}}$.

When the attacker does not know the parameters of underlying system such as network topology, line flow limits, each generator's capacity and ramp constraints, it is not possible for the attacker to find the optimal attacks, and the attacker just randomly chooses a set of loads to attack using either $\tilde{\underline{\mathbf{X}}}_j$ or $\tilde{\bar{\mathbf{X}}}_j$ for $j\in \tilde{\mathbb{L}}$. In the next section, we will show attacks on real-world load data upon IEEE standard systems that reveal the  vulnerabilities brought by either strategic or random load forecasting attacks.

 \begin{algorithm}
	\caption{Best-First Search}
	\label{algorithm}
	\begin{algorithmic}
		\REQUIRE Forecasted loads $L_1, L_2,...L_{|D|}\in \mathbb{R}^{|D|\times T}$
		\ENSURE Sets of Compromised Loads $\tilde{\mathbb{L} }\leftarrow \emptyset$
		\ENSURE Maximum Temperature perturbations $\epsilon$
		\ENSURE Number of compromised  loads $k=0$
		\STATE{Solve UC under normal forecasts }
		\WHILE{$G_t==\tilde{G}_t, t=1,...,T$ and $k\leq N_{adv}$}
		\STATE{Find the most vulnerable node $j$ and attack direction}
		\STATE{\texttt{Gradient Estimation} attacks to find $\tilde{L}_j$}
		\STATE{$\tilde{\mathbb{L} } \leftarrow j$}
		\STATE {Solve \eqref{equ:UC} to get $\tilde{G}_t$, $k+=1$}
		\ENDWHILE
		\STATE{Solve \eqref{equ:ED}} in real-time with $\tilde{G}_t$
	\end{algorithmic}
\end{algorithm}

%\subsection{Key Insights}
%For the general power system operations including a planning stage~(unit commitment) and a real-time operational stage~(economic dispatch), we observe the following characteristics of impacts by adversarial load forecasts:
%\begin{itemize}
%	\item Increasing the load maliciously will normally incur extra system costs, such as starting to operate redundant generators, using more expensive generation combinations and etc;
%	\item By decreasing the peak load maliciously, system operators would ignore the real peaks of future loads, and schedule fewer generators. This would potentially cause load shedding or failing to follow the severe ramps in the actual load patterns;
%	\item We assume an attacker with constrained capability on modifying the input features for load forecast models, and with no knowledge about the system parameters such as generator schedule and load forecasting model parameters. The proposed attack could be even more detrimental if the attacker possesses extra knowledge of the system and implement targeted attack during certain time periods.
%\end{itemize}

\section{Case Studies}
\label{sec:simulation}
%% !TEX root=main.tex
In this section, we show a detailed simulation on real-world Swiss load data, and show the threats posed  by our black-box data injection attacks on both the load forecasting algorithm itself and the power system operations. Thorough evaluations on 14-bus system indicate that proposed attacks cause load shedding with very high probability, while tests on 118-bus system indicate by compromising a small portion of nodal load forecasts, the attacker can cause security threats over the whole networks.

\subsection{Experimental Setup}
\textbf{Dataset Description}: We collected and queried hourly Swiss actual load data from European Network of Transmission System Operators for Electricity(ENTSO-E)'s API\footnote{https://transparency.entsoe.eu/} ranging from Jan 1st, 2015 to May 16th, 2017. The nominal load values are in the range of $[6,500MW,\; 9,500MW]$. We followed~\cite{marino2016building} to collect day-ahead historical weather forecasts coming from major cities in Switzerland such as Zurich, Basel, Lucerne and etc. All the weather data were queried from Dark Sky API~\footnote{https://darksky.net/forecast/47.3769,8.5414/us12/en}. We also collected other indicator features $\textbf{X}^{index}$, such as one-hot vectors of hour of day, weekdays and seasons. We evaluated the attack threats on split-out test data, and the forecast model parameters are kept away from the attacker throughout the black-box simulations. %We split $80\%$ of data as our training sets, and use the remaining $20\%$ of data on validating and evaluating the load forecasting prediction accuracy, attack performance and case studies on market operations. Note that even though we collected offline data to train and validate both of our load forecasting and attack models, these data collection procedures could be applied in an online fashion so that attacker could inject real-time adversarial attacks into certain load forecasting models.

\textbf{Power Systems Setup}:   We set up the IEEE 14-bus and IEEE 118-bus to study the system vulnerabilities brought by load forecasting attacks~\cite{christie2000power}. The grid has a total capacity of $15,500MW$ and $14,949.3MW$ respectively, which are both over $1.5$ times of yearly peak load.  We set the spinning reserve requirement as 3\% of the total forecasted demand based on~\cite{rebours2005spinning}. %During the run of day-ahead UC, either normal day-ahead forecasts or adversarial forecasts are used for generators scheduling; during the run of ED, the real loads are used for generation dispatch.
The models of UC and ED are implemented in Python using PyPSA~\cite{PyPSA}, and these two modules are directly interfaced with the load forecasting and attack algorithms using Tensorflow~\cite{abadi2016tensorflow}. No load shedding occurs when clean load forecasts are used.

\textbf{Model Training and Attack Implementation}: We set up an RNN with 3 layers for day-ahead load forecasting, and use standard stochastic gradient descent methods for model training~\cite{bottou2010large}.  Once the validation error converged, our RNN model reports an $1.58\%$ test error in mean absolute percentage error~(MAPE), which are comparable to the errors reported in several recent studies on load forecasting~\cite{kong2017short, quilumba2015using}. We use $L_\infty$ constraints on the attacker's capability \eqref{equ:constraint}, such that the attacker is constrained by the maximum deviation of perturbed temperature values. %We refer to Appendix for more training details and results on two other load forecasting models. %Based on discussion in Section~\ref{sec:market}, such accurate load forecasting performance is reliable to use in day-ahead unit commitment.
%We save the model parameters and keep them away from black-box attackers.
\begin{figure}[ht]
	\centering
	\includegraphics[width=1.0 \columnwidth]{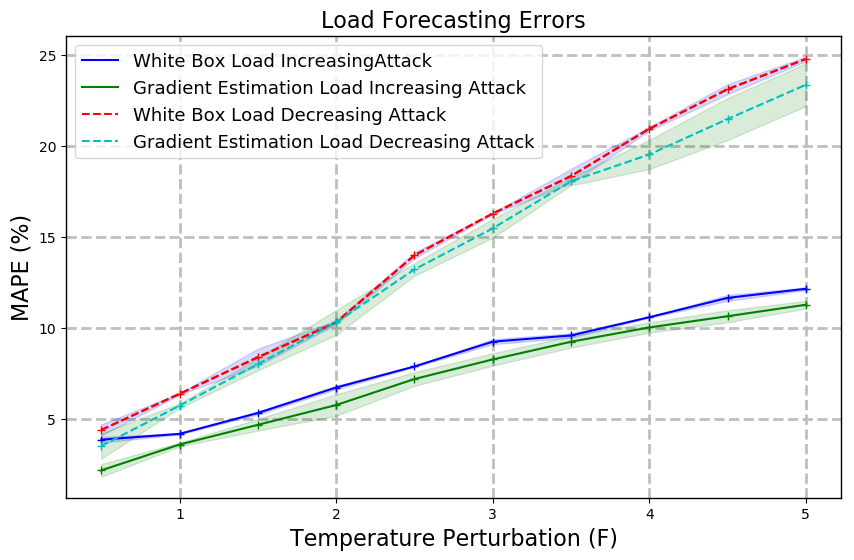}
	\caption{\small The forecast MAPE under attacks to increase the load and attacks to decrease the load using white box attacks and \texttt{gradient estimation} attacks. Simulation are run for three times with different random seeds, and shaded area denotes the variance.}
	\label{fig:errors}
\end{figure}

\subsection{Load Forecasting Performance}
We calibrate and compare the load forecasting model performance with and without adversarial attacks on test datasets. Though the forecasting model exhibits good performances on clean test data, we inject different level of perturbations generated by  \texttt{gradient estimation} method, and found the forecasting performance decrease drastically as the adversarial perturbations become larger. In Figure~\ref{fig:forecast_example} we visualize the RNN's load forecasting results for $300$ hours using \texttt{gradient estimation} algorithm with maximum perturbation on temperature of $1F$ and $5F$ respectively. The attacker tries to increase the load in the first $150$ hours, and to decrease the load in the latter hours. We observe that the algorithm finds the correct attack direction to either increase or decrease the load. What's more, with only $1F$ deviation on temperatures, the load forecasts changes over $500$MW at some time steps. When the attacker increases the perturbation to $5F$, large forecasts error over $1,200$MW are observed for the Swiss load. The temperature profile before and after attack still looks similar, which could avoid system operators' security inspection~(Figure \ref{fig:forecast_example}(b)). In Figure~\ref{fig:errors} we show by either attacking to increase the forecasted loads or attacking to decrease the forecasted loads, the \texttt{gradient estimation} attacks could achieve similar results compared to its white-box counterparts. More essentially, few degrees of malicious perturbations have caused large deviations on forecasted loads.

\begin{figure}[h]
	\centering
	\includegraphics[width=0.9 \columnwidth]{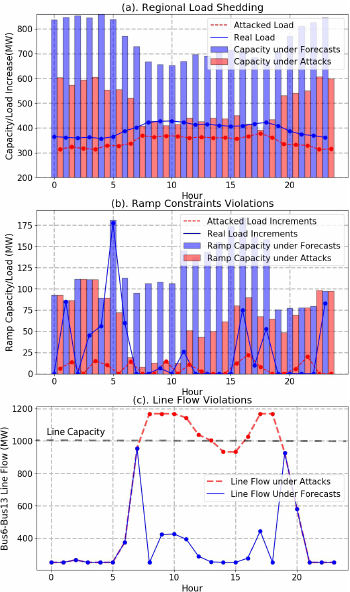}
	\caption{\small IEEE 14-bus threats visualization. (a). An example showing that forecasts under attack would cause generation limits violations when real loads exceed total generation capacity; bars indicate generators' available capacity. (b). An example showing that forecasts under attack would cause violation on ramp constraints during economic dispatch; bars indicate generators' available total up-ramp capabilities. Maximum allowed perturbations are $4F$. (c). An example showing that line flow exceeds limits.}
	\label{fig:UC}
\end{figure}

\begin{figure}[h]
	\centering
	\includegraphics[width=0.9 \columnwidth]{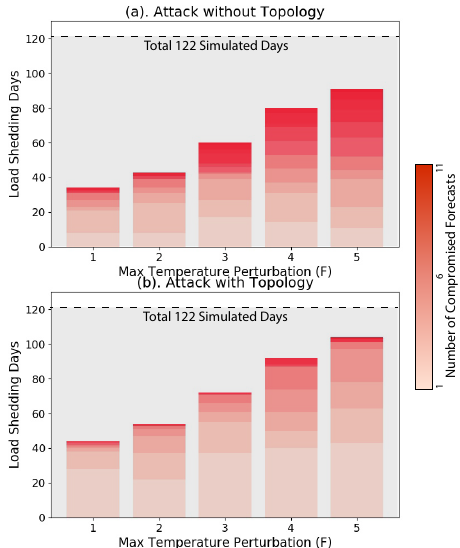}
	\caption{\small Load shedding caused by load forecasting attacks on IEEE 14-bus test case. For the 122 test days using Swiss load data, we compare the attack results when attacker either knows or does not know network topology. Larger perturbations on temperatures and extra information on topology helps the attacker to cause load shedding with fewer compromised loads.}
	\label{fig:min_load}
\end{figure}

\begin{figure*}[h]
	\centering
	\includegraphics[width=1.65 \columnwidth]{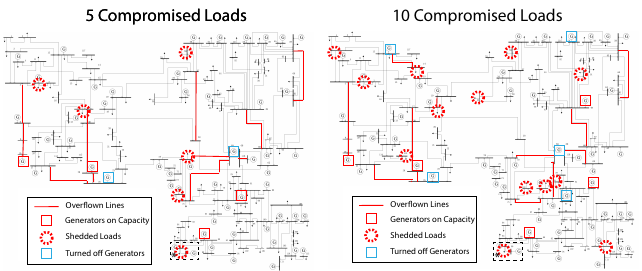}
	\caption{\small Threats posed by black-box attack on load forecasting algorithm in the IEEE 118-bus test case. By greedily attacking the nodal day-ahead load forecasts that leads to malicious commitment schedule, the attacker can incur a series of system operations threats such as load shedding, a number of overflown lines and on-capacity generators. For the bottom star-marked load shedded bus, the generator is off in both sub figures.}
	\label{fig:118}
\end{figure*}

\subsection{Impacts of Attacks on System Operations}
We find small adversarial perturbations over load forecasting input features even cause severe consequences on the power system operations. We are particularly interested in the case when $\tilde{G}_t$ is different from $G_t$ that causes infeasible solution for ED without load shedding. In Figure \ref{fig:UC}, we visualize different constraint violations when ED is solved.

We ran a thorough evaluation for $122$ test days on the 14-bus system, and evaluate if the attackers could cause load shedding under different level of knowledge and capabilities. In Figure \ref{fig:min_load}, we compare the number of days attacker cause load shedding with or without system topology information. In both cases, with larger perturbations added to the temperature forecasts, it is more likely to get an infeasible commitment schedule that causes load shedding. At one extreme, when the attacker knows network topology and is able to inject a $5F$ perturbation strategically on compromised load based on Algorithm~\ref{algorithm}, load shedding occurs on more than $100$ of the test days. At the other extreme, the attacker causes over $40$ days having shedded loads by only compromising one nodal load forecasts (Figure \ref{fig:min_load}(b)). More surprisingly, when the attacker does not know network topology and just selects compromised load randomly using either $\tilde{\underline{\mathbf{X}}}_j$ or $\tilde{\bar{\mathbf{X}}}_j$, the system operation is still very vulnerable to the proposed load forecasting attacks (Figure \ref{fig:min_load}(a)).

To furthur evaluate the attack's threats brought upon power system operation, one-day attack example on 118-bus test case is illustrated in Figure \ref{fig:118}. We assume the attacker is able to design greedy attacks based on known topology information in this case, and show that by only compromising a small subset of nodal load forecasts, the resulting day-ahead commitment schedule under attacks is shutting off several generators compared to that under normal forecasts. The network already observes a series of operation threats such as overflown lines, on-capacity generators and load shedding.

\section{Discussion and Conclusion}
\label{sec:conclusion}
In this paper, we studied the potential vulnerabilities generally existing in many load forecasting algorithms. Such vulnerabilities have been overlooked by the development of most if not all forecasting techniques. We designed a data injection attack which does not require parameters of the forecast algorithms, but leads to large increase in forecast errors. The proposed attack could adversarially  impact the decision making process for system operators. Experiments on real-world load datasets demonstrate such threats over power system operations. Such threats model along with damage analysis indicate that there needs more security evaluations in the design and implementation of load forecasting algorithms. In order to mitigate the damages brought by such false data injection attacks, countermeasures such as anomaly detection as well as other robust statistics are strongly recommended.

%\section*{Acknowledgement}
%The authors thank Baosen Zhang, Hossein Hosseini and Baicen Xiao from University of Washington, Shaohui Liu from Los Alamos National Laboratory for helpful discussion and valuable comments.

\bibliographystyle{IEEEtran}
\bibliography{ref}

\appendix
%% !TEX root=main.tex
\subsection{Details on \texttt{Learn and Attack} Algorithm}
\begin{figure}[h]
	\centering
	\includegraphics[width=1.0 \columnwidth]{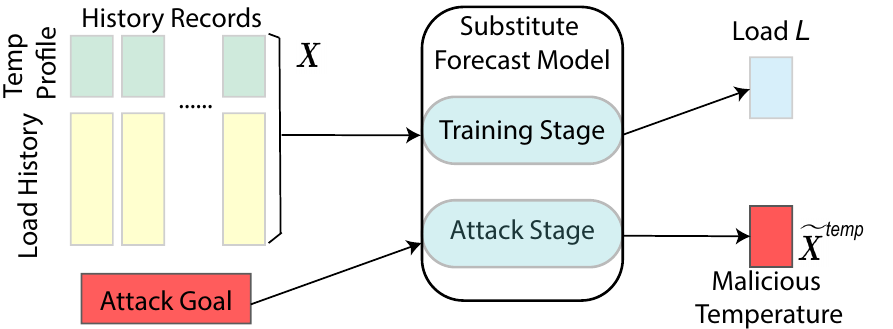}
	\caption{The attack schematic for \texttt{learn and attack} approach. During the training stage, the attacker uses historical data to learn a substitute forecast model; during attack stage, the attacker finds attack vector using substitute model and transfers it to unknown targeted system.}
	\label{fig:loadforecasting}
\end{figure}

In addition to the \texttt{gradient estimation} attack introduced in the main text, we also consider the \texttt{learn and attack} setting,  where we assume the attacker does not have access to the model parameters, and there is no query access to the model. The only knowledge the attacker has is a historical dataset $\tilde{\mathcal{D}}_{tr}$, which includes same features under same distributions as data set  $\mathcal{D}_{tr}$ used to train the load forecasting model\footnote{In Learn and Attack setting, we make assumption that the attacker know the family of targeted load forecasting model,  e.g., a feedforward neural networks or a Recurrent Neural Networks.}. The proposed attack algorithm consists of a \emph{training phase} and an \emph{attack phase} as shown in Figure~\ref{fig:loadforecasting}. In the training phase, the attacker trains substitute model $f_{\tilde{\theta}}$ based on $\tilde{\mathcal{D}}_{tr}$ to minimize the training loss.  In the attack phase, the attacker pretends that the substitute model is the true load forecast model and performs white-box attacks on it to find the attack vectors. This strategy is based on the assumption that the substitute model behaves similarly to the true model not only for the training data $\textbf{X}$, but also for the attack vector $\tilde{\textbf{X}}$. Then by injecting $\tilde{\textbf{X}}$ into the true load forecasting model, the forecast values go to attacker's desired directions.

\begin{figure*}[h]
	\centering
	\includegraphics[width=2.0 \columnwidth]{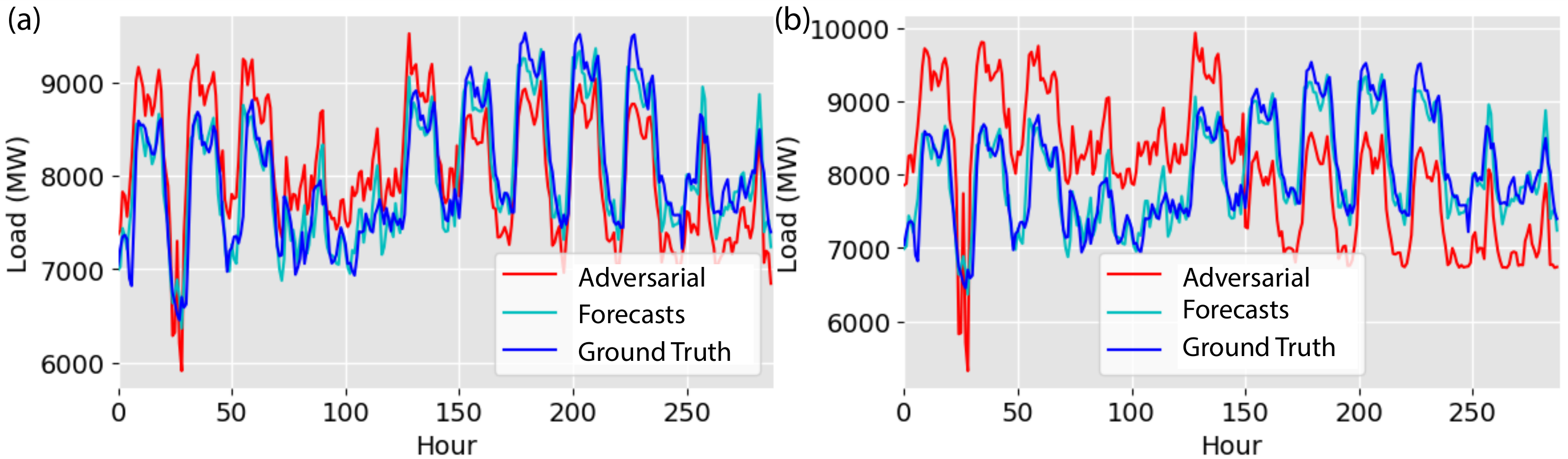}
	\caption{We show 300 hours forecasts based on original and adversarial temperature data for the aggregated load of Switzerland. The load forecasting algorithm is an recurrent neural networks with inputs composed of past load, regional temperature
forecast values and weather indicators. The attack perturbations are generated by using the \texttt{learn and attack} method, and
it implements load maximization strategy in the first 150 hours and load minimization strategy in the latter 150 hours. (a).
Load forecasting results with temperature attack constraint of (maximum perturbations) 1F ; (b). load forecasting results with
temperature attack constraint of 5F .}
	\label{fig:forecast_example2}
\end{figure*}

It is useful to evaluate the \emph{transferability} of proposed attacks across different set of models $f_\theta$ and $f_{\tilde{\theta}}$. The phenomenon of transferability in adversarial attacks for machine learning models have been discussed in  \cite{papernot2016transferability, hosseini2017blocking}, where adversarial instance generated using $f_{\tilde{\theta}}$ can be also treated as an adversarial instance by $f_{\theta}$ with high probability. The theoretical understanding of why attacks transfer remains an open question and is out of scope for this paper. In Figure~\ref{fig:forecast_example2} we show such \emph{transferability} also exists in the load forecasting model using same test case on Switzerland load forecasting. The temperature inputs are generated by implementing the iterative gradient update based on a substitute model under $L_\infty$-norm of attack perturbations, yet such adversarial temperature values also mislead the (unknown) true load forecasting model to be wildly inaccurate.

\subsection{Details on Best First Search Attack Selection}
As described in Algorithm~\ref{algorithm} in Section \ref{subsec:strategy}, since it is computationally inefficient and challenging for the attacker to solve a tri-level attack problem (attacker solves an adversarial UC \eqref{equ:UC} using $\tilde{L}_q^t$, operator minimizes ED costs \eqref{equ:ED} using $\tilde{G}_t$, attacker maximizes ED costs \eqref{equ:ED} using $\tilde{G}_t$), we propose a modified best-first search algorithm~\cite{dechter1985generalized} for attackers to compromise a limited number of nodal forecasts. Essentially, in order to solve the tri-level problem expressively to find attack vectors, our proposed attack on networks iteratively finds the most vulnerable node to inject attacks.

Suppose for $|Q|$ loads in the network, the attacker can at most $N_{adv}$ nodal load forecasts to avoid detection by system operators. Meanwhile, the constraints on attacker's capabilities described in Section \ref{sec:method} shall hold throughout data perturbation attacks. If an attacker could cause the day-ahead UC schedule $\tilde{G}_t, t=1,...,T$ using $\tilde L_q^t, q=1,...,|Q|,t=1,...,T$ shift from $G_t, t=1,...,T$, then it is expected the solution of ED will change. There are several possible circumstances by using $\tilde{L}_q^t$:
\begin{itemize}
	\item Increasing the load maliciously will possibly incur extra system costs, such as starting to operate redundant generators, using more expensive generation combinations and etc;
	\item Decreasing the load maliciously will possibly incur infeasble generation schedules during real-time dispatch, since there may be fewer generators scheduled than normal conditions, which will cause generators reach capacity or line flow exceed limits;
	\item Decreasing the peak value of load maliciously will possibly cause UC solver ignore peak values, which may cause generators reach ramping limits.
\end{itemize}

In our proposed algorithm for attackers to find most vulnerable nodal forecasts and inject attack forecasts, we design an iterative search scheme. At each iteration, the attacker solves day-ahead UC \eqref{equ:UC} based on current $\tilde{L}$, and check which generator's schedule is most prone to change. Then the attacker decides the next compromised load node $j$. Since it is more possible to change UC schedule by chaning the load profile to greater extents, the attacker then proposes two attack samples $\tilde{\underline{\mathbf{X}}}_j$ and $\tilde{\bar{\mathbf{X}}}_j$ which minimizes and maximizes the nodal forecasts respectively. Such adversarial nodal load forecast is inserted into next iteration's UC problem. 

In our implementation, we simply find the prone-to-change generation schedule based on following criteria:
\begin{itemize}
	\item The generator's dispatch during day-ahead UC is either approaching to generation limit or reaching zero output;
	\item The line flow during day-ahead UC is approaching line flow limits;
\end{itemize}

And simulation results indicate such search algorithm is quite efficient to find the attack injections. Other criteria such as ramping constraints and generation on-off limits could be also utilized to find the most vulnerable nodal forecasts.

For the case when system topology and parameters are unknown, the attacker just selects a random set of compromised load and inject attacks perturbations using either $\tilde{\underline{\mathbf{X}}}_j$ or $\tilde{\bar{\mathbf{X}}}_j$ for the compromised nodal loads. The result shown in Fig. \ref{fig:min_load} validates that in the case without topology knowledge, compromising the same number of loads are causing fewer load shedding days than in the case when topology is known.

\subsection{Details on Attacks Implementation}
\label{sec:app}
In addition to the results we have shown in the main texts using RNN as load forecasting algorithm, we added ablation study on different load forecasting algorithms, as well as attack performance and computation time analysis.

\subsubsection{Forecasting Method}

\textbf{Feed-Forward Neural Networks}
A multi-layered, feed-forward neural networks~(NN) has been widely used to represent the nonlinearities between input features and output forecasts~\cite{hippert2001neural}. For the input layer of neural networks, each neuron represents one feature of training input, and all features of past $H$ steps $(\textbf{X}_{t-H},...,\textbf{X}_t)$ are stacked as the inputs. For each intermediate layer, NN could have a tunable number of hidden units, which represent the input feature combinations.

\noindent \textbf{Recurrent Neural Networks}
As described in main texts, RNN feeds each step's input $\textbf{X}_t$ sequentially, and outputs a hidden unit to represent the feature combination of current input and historical features. The last neuron outputs the forecasted load values in the load forecasting attack~\cite{chen2018machine}.

\begin{figure}[!h]
	\centering
	\includegraphics[width=0.95 \columnwidth]{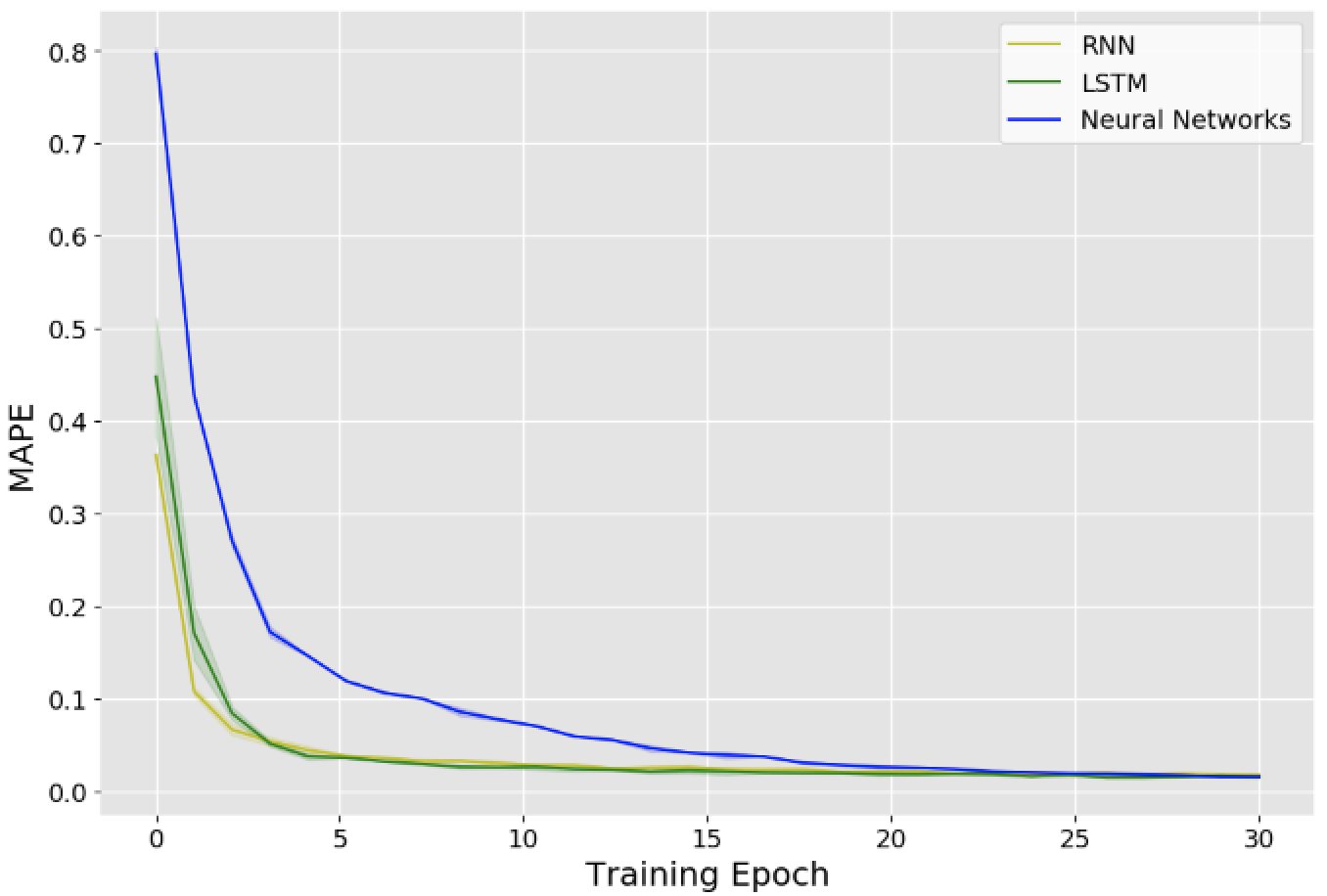}
	\caption{All three forecasting models, show convergence of forecast error on validation data as  training evolves. Shaded areas show the variance of MAPE.}
	\label{fig:forecast_training}
\end{figure}

\noindent \textbf{Long Short-Term Memory}
Long Short-Term Memory network~(LSTM) is designed to deal with the vanishing gradient problem existing in the RNN with long-time dependencies~\cite{kong2017short,chen2018short}. The major improvements over RNN are the design of ``forget" gates to model the temporal dependencies and capture long time dependencies in load patterns more accurately.

\begin{table}[h]
	\centering
	% \begin{tabular}{>{\centering}m{5.8cm} >{\centering}m{1.2cm} >{\centering}m{1.2cm} >{\centering\arraybackslash}m{1.2cm}}
	\begin{tabular}{cccc}
		\toprule[0.4mm]
		Forecasts Models&NN&RNN&LSTM\\
		\toprule[0.4mm]
		%\hline
		Number of Layers & 4 & 3 & 3\\
		%\hline
		Training Epochs & 20 &30 &30 \\
		%\hline
		Hidden Units in First Layer& 512& 64&64\\
		\bottomrule[0.4mm]
	\end{tabular}
	\caption{Model architectures and training configurations for load forecasting algorithms used in the simulations.}
	\label{table:nn}
\end{table}

\begin{figure*}
	\centering
	\includegraphics[width=2.0\columnwidth]{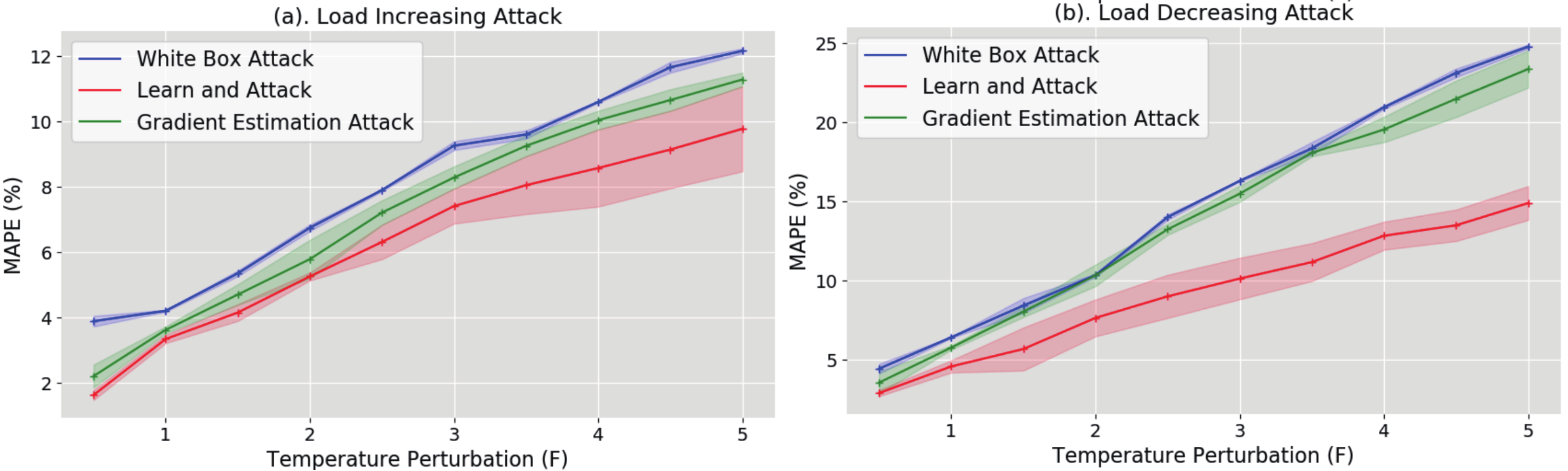}
	\caption{The forecast MAPE under (a). attacks to increase the load; and (b). attacks to decrease the load. Simulations are run for three times with different random seeds for all three attacks, and
		shaded area denotes the variance.}
	\label{fig:accuracy}
\end{figure*}

\subsubsection{Training and Attack Details}
We set up all load forecasting models using Tensorflow~\cite{abadi2016tensorflow} package in Python. Standard model architectures such as Dropout layers and nonlinear activation functions (e.g., ReLU or Sigmoid functions) are adopted in the deep learning models~\cite{nair2010rectified}. Since all three networks are set up to solve the load forecasting regression problem, we set the first layer having most neurons, and decrease the number of units in subsequent layers. Table~\ref{table:nn} records our load forecasting model setup. We split $80\%$ of data as our training sets, and
use the remaining $20\%$ of data on validating and evaluating the
load forecasting prediction accuracy, attack performance and case
studies on market operations.  Such data collection procedures could also be applied in
an online fashion so that attacker could inject real-time adversarial
attacks into certain load forecasting models.

\begin{table}
	\centering
	\begin{tabular}
		{P{1.8cm}| P{1.8cm}P{1.8cm}P{1.8cm}}
		Forecasts Error (MAPE)&Clean Data&Learn and Attack&Gradient Estimation\\
		\toprule[0.4mm]
		%\hline
		NN & $1.68\%$ & $12.72\%$ & $13.09\%$\\
		%\hline
		RNN & $1.58\%$ &$9.82\%$ &$11.68\%$ \\
		%\hline
		LSTM & $1.51\%$& $9.04\%$&$11.87\%$\\
		\bottomrule[0.4mm]
	\end{tabular}
	\caption{Forecasts errors evaluated on clean test data and adversarial data for 3 different forecast models. Allowed maximum perturbations are $4F$.}
	\label{table:forecast_error}
\end{table}

As shown in Figure \ref{fig:forecast_training}, all three load forecasting algorithms' validation loss are converged during training, and we use the trained model in the subsequent planning and operation problem as well as the testbed for attack algorithms. Plots are showing the mean and variance during 3 runs.

For the substitute model training of \texttt{learn and attack} method, we keep the training set $\tilde{\mathcal{D}}$ same as the load forecasting model training set $\mathcal{D}$. Decreasing the size of $\tilde{\mathcal{D}}$ or using different substitute dataset could decrease the performance of \texttt{learn and attack}. Table \ref{table:forecast_error} compares all three load forecasting models' performance using clean and adversarial data. For both \texttt{learn and attack} and \texttt{gradient estimation} algorithms, they distort all three load forecasting models' output and increase model's forecast error. \texttt{Gradient estimation} attack works generally better for all three models, and this is due to estimating the gradients via querying $f_\theta$ directly is more accurate than calculating it from the substitute model and transferring to $f_\theta$.

In Figure~\ref{fig:accuracy}, we evaluate RNN's load forecasting performance under two attack strategies: load maximization or load minimization. We observe \texttt{gradient estimation} attack causes similar MAPE compared to white box attack.  The load decreasing attack is normally more successful than load increasing attack in terms of MAPE. Load minimization attack is more harmful results than load increasing ones, since increased forecasts only let system operators start up more generations, while adversarially decreasing the forecasted load leads to wrong generation decisions that fails to meet the larger real load.

\begin{table}[!h]
	\centering
	% \begin{tabular}{>{\centering}m{6.2cm} >{\centering}m{1.2cm} >{\centering}m{1.2cm} >{\centering\arraybackslash}m{1.2cm}}
	\begin{tabular}{cccc}
		Forecasts Models&NN&RNN&LSTM\\
		\toprule[0.4mm]
		%\hline
		Training Time & 12.988 & 47.998 & 143.830\\
		%\hline
		\texttt{Learn and Attack} & 0.133 &0.157 &0.579 \\
		%\hline
		\texttt{Gradient Estimation Attack} & 0.082 & 0.119 &0.253\\
		\bottomrule[0.4mm]
	\end{tabular}
	\caption{Computation time (in seconds) for load forecasting model training and implementation time for attacks.}
	\label{table:time}
\end{table}

\subsubsection{Computation Time}
\label{sec:app2}
We recorded the computation time for neural network training and the implementation time for two proposed attack algorithms. All time are recorded on a laptop with Intel 2.3GHz Core i5-8259U 4 Cores CPU and 8 GB RAM. The training time for NN, RNN and LSTM are calculated for $ 20, 30$ and $30$ epochs respectively. The implementation time for the attacks are averaged over all test instances. We observed that \texttt{learn and attack} approach takes longer time than  \texttt{gradient estimation} due to the longer time taken to calculate gradient signs over the whole neural networks; and as LSTM includes more complicated model architectures, it takes longer time to find the adversarial instance. Yet compared to the long model training time and application scenarios of day-ahead forecasts, the attacker is still efficient enough to find the adversarial perturbations. Such efficient computation enables attacker to find the most vulnerable loads or to attack the system operations in very short time.

\subsubsection{Code Availability}
The implementation code for forecasts, attacks, and power system operations are all available at \url{https://github.com/chennnnnyize/load_forecasts_attack}.

\end{document}